\documentclass[reqno]{svproc}
\usepackage{amssymb,amsmath}
\usepackage{graphicx}
\usepackage[left]{lineno}

\usepackage{algorithm}
\usepackage{algpseudocode}
\usepackage{booktabs}
\usepackage{longtable}
\usepackage{caption}
\usepackage{soul}
\usepackage{cite}
\usepackage{float}

\DeclareMathOperator*{\argmax}{arg\,max}

\newcommand{\stepsize}[0]{h}

\usepackage[margin=2.4cm]{geometry}

\usepackage{url}

\begin{document}
	\mainmatter              
	\title{Feature weighting for data analysis via evolutionary simulation}
	\titlerunning{}  
	%
	\author{Aris Daniilidis\inst{1} \and Alberto Dom\'inguez Corella\inst{1,2} \and Philipp Wissgott\inst{3}}
	
	\authorrunning{Aris Daniilidis, Alberto Dom\'inguez Corella, and Philipp Wissgott}
	
	\tocauthor{Aris Daniilidis, Alberto Dom\'inguez Corella, and Philipp Wissgott}
	
	\institute{
		Institut f\"{u}r Stochastik und Wirtschaftsmathematik, Variational Analysis, Dynamics and Operations Research Unit E105-04, Technische Universit\"at Wien, Wiedner Hauptstra\ss{e} 8, 1040 Vienna, Austria\\
		\email{aris.daniilidis@tuwien.ac.at}\\
		\email{alberto.of.sonora@gmail.com}
		\and
		Institut f\"{u}r Mathematik und Wissenschaftliches Rechnen, Universit\"{a}t Graz, Heinrichstra{\ss}e 36, A-8010 Graz, Austria
		\and
		danube.ai solutions gmbh,\\
		1040 Vienna, Austria\\
		\email{philipp@danube.ai}
	}

	\maketitle              
	
	\begin{abstract}   
We analyze an algorithm for assigning weights prior to scalarization in discrete multi-objective problems arising from data analysis. The algorithm evolves weights (interpreted as the relevance of features) by a replicator-type dynamic on the standard simplex, with update indices computed from a normalized data matrix. We prove that the resulting sequence converges globally to a unique interior equilibrium, yielding non-degenerate limiting weights. 
		\keywords{ feature weighting, evolutionary algorithms, multi-objective optimization, data  analysis.}
	\end{abstract}

	\section{Introduction, data problem, and theoretical results}
   In many applied domains—such as economics, energy systems, technological design, and consumer recommendation—the relative importance of features describing alternatives is often unknown and must be inferred from data or performance criteria.
   \smallbreak\noindent
   We analyze an algorithm recently introduced in \cite[Section III]{W_2025} for estimating the relative importance of dataset features. The feature weights are modeled as probability vectors and are iteratively updated from an initial distribution according to a replicator-type dynamic. The method is inspired by evolutionary game theory: features are interpreted as players whose weights evolve according to their relative performance, whereas the update rules depend on strategy functions determined by the data type. We restrict attention to two simple strategies: in the first strategy, features with strong variation should increase in weight; in the second one, datasets penalize over-reliance on a single feature. 
    \smallbreak\noindent
    We next describe the feature-weighting problem, present the algorithm and its convergence results, and defer proofs to Section~\ref{proofs}. A numerical example appears in Section~\ref{numexp}, and connections to genetic algorithms and evolutionary systems are discussed in Section~\ref{sec:Analogies}.
	\subsection{The data problem}\label{sdp}
	Consider a finite set $\mathcal X=\{x_1,\dots, x_n\}\subseteq \mathbb R^m$ representing the available choices of a rational agent. The set can also represent different datasets. For each $j\in\{1,\dots,m\}$, the number $x_{ij}\in\mathbb R$ represents a feature of option/dataset $i\in\{1,\dots,n\}$. The agent seeks to choose an option that maximizes some features (benefit-related) and minimizes others (cost-related), however this is not always possible and one seeks a combined concept of \textit{optimality} in this context.  This problem falls within the class of so-called multi-objective problems. While this framework offers several notions of solution—such as Pareto optima—these do not capture the relative importance of each feature. In this paper, we analyze the possible relevance (quantified as an index) of features  that can be obtained from the  data of the problem. 
	\smallbreak\noindent 
	Below, we consider an algorithm that outputs a probability vector $\gamma \in \mathbb{R}^m$, where  each $\gamma_j$ represents an index of the relevance in optimizing feature $j \in \{1, \dots, m\}$. The algorithm is inspired by evolutionary game theory, drawing an analogy in which datasets correspond to \textit{organisms} and features to \textit{genes}; this interpretation is described in detail in Section \ref{sec:Analogies}. 
	\smallbreak\noindent 
	The input of the algorithm can be thought of as a matrix
	\begin{align*}
		X:= \begin{bmatrix}
			x_{11} & x_{12} & \dots & x_{1m} \\
			x_{21} & x_{22} & \dots & x_{2m} \\
			\vdots & \vdots & \ddots & \vdots \\
			x_{n1} & x_{n2} & \dots & x_{nm}
		\end{bmatrix}.
	\end{align*}
	To effectively process data and provide a consistent interpretation across features, the raw data is transformed into a normalized representation. The purposes of this normalization are twofold.  First, we want to ensure that all feature values lie in a consistent scale, and second, we want to provide data with a percentage-based interpretation, thereby requiring $\Phi_{ij} \in [0,1]$ for all $i \in \{1, \dots, n\}$ and $j \in \{1, \dots, m\}$. Once the columnwise normalizations are fixed for a problem class, any dataset with the same structure can be processed dynamically without further manual tuning.
	\smallbreak\noindent 
	Consider the normalization represented by a matrix 
	\begin{align*}
		\Phi:= \begin{bmatrix}
			\Phi_{11} & \Phi_{12} & \dots & \Phi_{1m} \\
			\Phi_{21} & \Phi_{22} & \dots & \Phi_{2m} \\
			\vdots & \vdots & \ddots & \vdots \\
			\Phi_{n1} & \Phi_{n2} & \dots & \Phi_{nm}
		\end{bmatrix}
	\end{align*}
	with entries in $[0,1]$. Such normalization transforms the  optimization of features to a multi-objective problem. 
    \begin{example}
        When dealing with matrices $X$ containing non-negative data, for a column $j\in\{1,\dots,m\}$ not identically zero, a possible transformation applied to a cost-related feature could be given by
	\begin{align*}
		\Phi_{ij}\,=\,1\,-\,\frac{x_{ij}}{\displaystyle\max_{s\in\{1,\dots,n\}} x_{sj}}\,,\qquad \text{for all \,}\, i\in\{1,\dots,n\}. 
	\end{align*}
	This transformation applies an \textit{inverted  normalization} to the column. Each entry is scaled by its feature's maximum, and this ratio is subtracted from $1$. Consequently, $\Phi_{ij} \in [0, 1]$, where $0$ denotes the feature's maximum value (no deviation) and $1$ signifies an entry of $0$ (100\% deviation). This yields a measure that is  useful when the objective is to minimize a feature rather than maximize it. 
    \end{example}
    \noindent
   We emphasize that the normalization is specified separately for each column \(j\in\{1,\dots,m\}\), and then kept fixed across all \(i\in\{1,\dots,n\}\). Choosing an appropriate normalization for column \(j\) requires an understanding of the feature encoded by that column. In the evolutionary interpretation introduced in Section~\ref{sec:Analogies}, this corresponds to determining the local fitness of feature \(j\). Qualitatively, one must decide whether higher values of that feature are favorable or unfavorable for solutions \(x_i\).
       \smallbreak\noindent
    Following the so-called \textit{No Free Lunch} theorem—here interpreted as a metatheorem stating that no optimization algorithm is universally superior or able to solve all problems—we remark that the relevance of the algorithm  lies in its behavior on structured datasets, where inductive biases reflecting feature interactions can be effectively exploited. 
	
	\subsection{Algorithm and convergence results}\label{Althm} 
	Let us begin with some notation. 
	We denote by  \[
	\mathcal{K}^{m-1} := \left\{ \gamma \in \mathbb{R}^m \Big| \gamma_j \geq 0 \,\,\, \forall j\in\{1,\dots, m\} \,\,\,\text{and}\,\,\, \sum_{j=1}^m \gamma_j = 1 \right\}
	\]
	the standard simplex in \( \mathbb{R}^m \). The elements of $\mathcal K^{m-1}$ will correspond to the scalarization vectors (feature weights) of the multi-objective  problem. We write
	\begin{align}\label{eq:meandef}
		\widetilde{\Phi_j}:=\frac{1}{n}\sum_{i=1}^n\Phi_{ij}
	\end{align}
	for the column average of the matrix $\Phi \in M_{n \times m}(\mathbb{R})$. Notice that each average lies in the interval $[0,1]$. 
    \smallbreak\noindent
	For each $ j\in\{1,\dots, m\}$,  consider the functions $\Delta_j^{\operatorname*{dom}},\Delta_j^{\operatorname*{bal}}:\mathcal K^{m-1}\to\mathbb R$ given by 
	\begin{align}\label{eq:principledeltadef}
		\Delta_j^{\operatorname*{dom}}(\gamma):=\gamma_j\Big(\widetilde{\Phi_j}-\frac{1}{2}\Big) \quad \text{and}\quad \Delta_j^{\operatorname*{bal}}(\gamma):= -2\Big(\gamma_j\widetilde{\Phi_j} - \frac{1}{m}\sum_{s=1}^m\gamma_s\widetilde{\Phi_s}\Big).
	\end{align}
	The quantities $\Delta_j^{\operatorname*{dom}}$ (dominance term) and $\Delta_j^{\operatorname*{bal}}$ (balanced term) can be viewed as indices that evaluate the contribution of feature $j\in\{1,\dots,m\}$  under different criteria. The dominance term measures whether a feature  tends to be high or low on average. It is positive when the column average is above $0.5$, and its effect grows with the weight given to that feature. In this sense, it rewards features that are both strongly weighted and above the midpoint of the scale.
	The balance term compares the weighted contribution of a feature  with the overall weighted mean across all features. It is positive when the feature contributes less than the global mean, and negative when it contributes more. This acts as a correction; it discourages features from dominating too strongly and favors those that are closer to balance with the rest
	\smallbreak\noindent 
	We also consider the function  $\Delta_j:\mathcal K^{m-1}\to\mathbb R$ given by
	\begin{equation}\label{eq:tri}
		\Delta_j(\gamma):=\Delta_j^{\operatorname*{dom}}(\gamma)+ 	\Delta_j^{\operatorname*{bal}}(\gamma) = -\gamma_j\Bigl(\widetilde{\Phi}_j+\tfrac12\Bigr)
	\;+\;\frac{2}{m}\sum_{s=1}^m\gamma_s\widetilde{\Phi_s}\,.
	\end{equation}
	This is an index of the relative contribution of feature $j\in\{1,\dots,m\}$, combining the effects of the dominance and balance terms.
	\smallbreak\noindent
    In what follows, we assume that there are at least two features under consideration, i.e., that $m\ge 2$. All along the paper, we fix the quantities 
    $$\widetilde{\Phi}_{\operatorname{max}}:=\max_{j\in\{1,\dots,m\}}\widetilde{\Phi_j}$$ and 
    \begin{equation} \label{eq:h0}
           h_0:=\left(\frac12+\left(1-\frac{2}{m}\right)\widetilde{\Phi}_{\operatorname{max}}\right)^{-1}>0.
    \end{equation}
	We consider a fixed step size (learning rate) parameter $h \in (0, h_0)$. Given an initial weight \( \gamma^0 \in \mathcal{K}^{m-1} \), we say that $(\gamma^k)_{k\in\mathbb N}\subseteq \mathcal K^{m-1}$ is a sequence generated by Algorithm \ref{gAI} if, for each $k\in\mathbb N$,
	\begin{align}\label{eq:replicatorequation1}
		  \gamma_j^{k+1} = \frac{\gamma_j^k \left( 1 + h\Delta_j(\gamma^k) \right)}{\displaystyle\sum_{s=1}^m \gamma_s^k \left( 1 + h\Delta_s(\gamma^k) \right)} \quad \forall j\in \{1, \dots, m\}.
	\end{align}

    \begin{algorithm}
		\caption{\label{gAI}}
		\begin{algorithmic}[1]
			\Require Initial weight \( \gamma^0 \in \mathcal{K}^{m-1} \).
			\State Initialize \( \gamma \longleftarrow \gamma^0 \)
			\For{$k \in\mathbb N$}
			\State Compute \( \Delta \longleftarrow\Delta(\gamma) \)
			\For{$j \gets 1$ to $m$}
			\State Update \( \gamma_j \longleftarrow \gamma_j (1 + h\Delta_j) \)
			\EndFor
			\State Normalize \( \gamma \longleftarrow \gamma / \sum_{s=1}^m \gamma_s \)
			\EndFor
		\end{algorithmic}
	\end{algorithm}
\noindent The update map in Algorithm~\ref{gAI} is a multiplicative renormalization on the simplex. 
\begin{remark}
The use of the step size $h$ in \eqref{eq:replicatorequation1} ensures that Algorithm~\ref{gAI} is well-defined. Should this damping term be omitted, the factor $1+\Delta_j(\gamma)$ might not be positive: indeed for $h=1$, take $m=5$,
\[
\widetilde\Phi=\left(1,\frac12,\frac12,\frac12,\frac12\right)\quad \text{and}\quad 
\gamma^0=\left(\frac{15}{16},\frac1{64},\frac1{64},\frac1{64},\frac1{64}\right).
\]
Then a simple computation shows that $1+\Delta_1(\gamma^0)=-\frac{3}{160}<0$ and the corresponding multiplicative update does not define a point of the simplex.
\end{remark}
The following lemma shows the well-posedness of the algorithm with the damping term.
\begin{lemma}[Well-posedness of the algorithm]\label{lem0}
For every $\gamma\in\mathcal K^{m-1}$, there holds
\[
\sum_{s=1}^m\gamma_s\bigl(1+h\Delta_s(\gamma)\bigr)>0
\quad\text{and}\quad
\gamma_j\bigl(1+h\Delta_j(\gamma)\bigr)\ge 0
\qquad \forall j\in\{1,\dots,m\}.
\]
\end{lemma}
    \noindent
	We see then that the update rule defines a discrete dynamical system on the simplex. At each step, the weight of feature $j\in\{1,\dots, m\}$ is multiplied by a factor proportional to $1+h\Delta_j(\gamma^k)$, and the whole vector is renormalized so that the updated weights again lie in $\mathcal{K}^{m-1}$. 
Features with positive indices $\Delta_j(\gamma^k)$ are reinforced, while those with negative indices are diminished. This mechanism is reminiscent of the spirit of replicator dynamics, where the evolution of a distribution is driven by relative performance, with better-performing features gaining weight at the expense of weaker ones.
\begin{proposition}[Non-degeneracy of limit points of the algorithm]\label{prop0}
		Consider the set
		\[
		\mathcal K_+^{m-1} := \left\{ \gamma \in \mathbb{R}^m \Big| \gamma_j > 0 \,\,\, \forall j\in\{1,\dots, m\} \,\,\,\text{and}\,\,\, \sum_{j=1}^m \gamma_j = 1 \right\}
		\]
		Let $(\gamma^k)_{k\in\mathbb N}$ be a sequence generated by Algorithm \ref{gAI}. If $\gamma^0\in\mathcal K^{m-1}_+$, then $(\gamma^k)_{k\in\mathbb N}$, as well as any accumulation point of this sequence, also belong to $\mathcal K_+^{m-1}$. 
	\end{proposition}
	 From the geometric point of view, the dynamics avoid collapsing onto the boundary of the simplex; thereby ensuring that no feature is ultimately discarded.  
	With this invariance established, we now turn to the asymptotic behavior of the sequence. The key question is whether the updates converge to a stable distribution inside the simplex, rather than oscillating or drifting indefinitely.
	\begin{theorem}[Convergence of the algorithm to an equilibrium]\label{thm0}
		 Set
        \[
        a_j:=(\widetilde{\Phi_j}+1/2)^{-1},\quad j\in \{1.\dots,m\}, \quad 
        A:=\sum_{j=1}^m a_j.  
        \]
        Let $\gamma^*\in\mathcal K_{+}^{m-1}$ be given by 
		\begin{align}\label{fp}
			\gamma_j^* = \displaystyle \frac{a_j}{A},\qquad \text{for all } \,j\in\{1,\dots,m\}. 
		\end{align}
        Consider the entropy function $V:\mathcal K_+^{m-1}\to [0,+\infty)$ given by
	\begin{align}\label{eq:Lyapunovfunction}
		V(\gamma):=\sum_{j=1}^m \gamma_j^*\log\Big(\frac{\gamma_j^*}{\gamma_j}\Big).
	\end{align}
    If $(\gamma^k)_{k\in\mathbb N}$ is a sequence generated by Algorithm \ref{gAI} satisfying $\gamma^0\in\mathcal K_+^{m-1}$, then the following properties hold. 
    \begin{itemize}
        \item[$(a)$] $V(\gamma^*)=0$ and $V(\gamma)>0$ for all $\gamma\in\mathcal K^{m-1}_+\setminus\{\gamma^*\}$;
        \smallbreak
         \item[$(b)$] $V(\gamma^{k+1})\le V(\gamma^{k})$ for all $k\in\mathbb N$, in particular, $V$ is a Lyapunov function of the algorithm;
        \smallbreak
         \item[$(c)$] $\gamma^k\longrightarrow\gamma^*$ as $k\longrightarrow+\infty$.  
    \end{itemize}
	\end{theorem}
     Formula \eqref{fp} reveals that the limiting weight of each feature depends monotonically on its column average; it decreases as the average increases. In particular, features with larger averages receive less weight, while those with smaller averages receive more. We note that the use of Lyapunov functions is a well known technique to analyze the dynamics in evolutionary games  \cite{Hofbauer1998}.
		\begin{figure}[h!]
		\centering
		\includegraphics[width=0.7\textwidth]{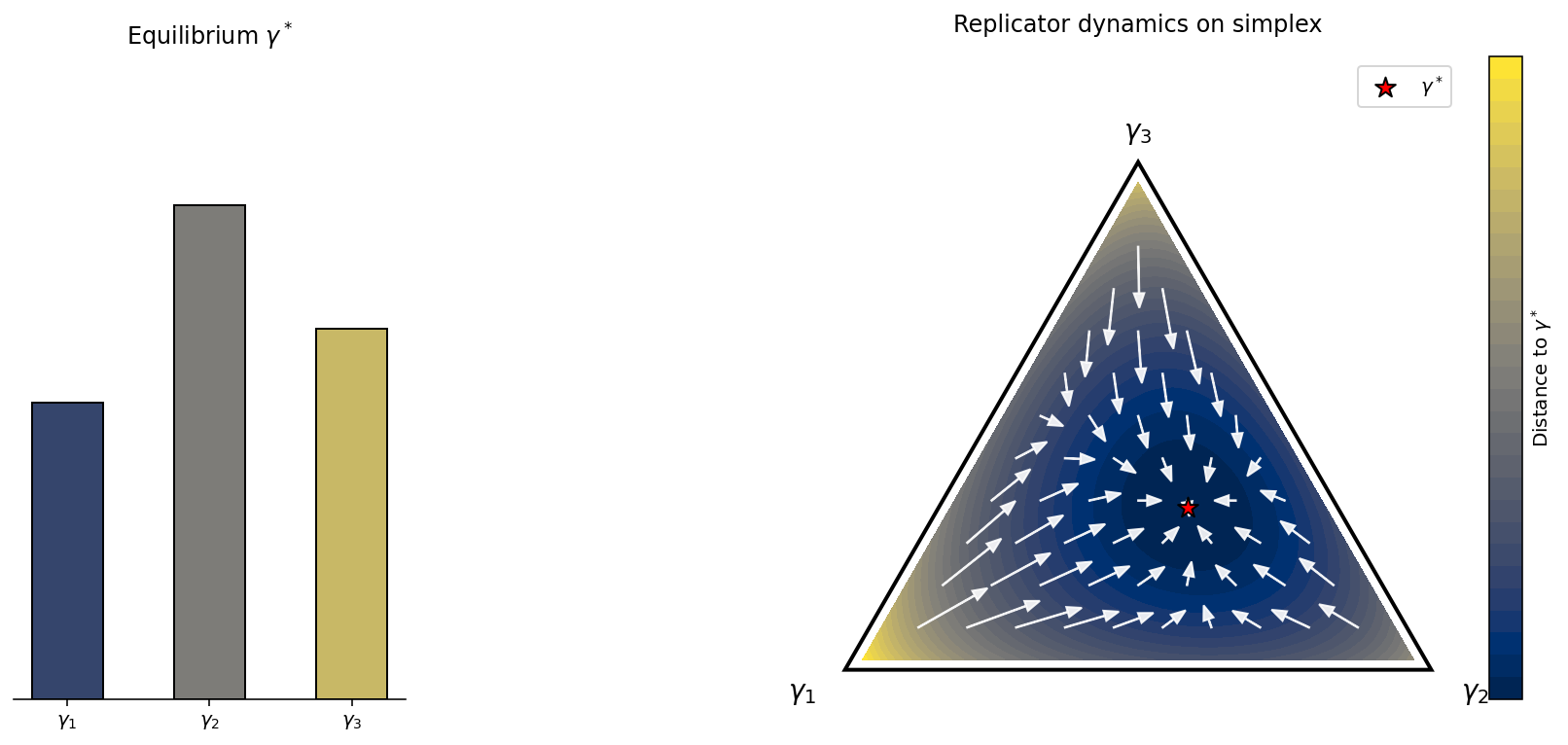}
		\caption{Illustration of the replicator-type feature weighting algorithm. 
		}
		\label{fig:triangle}
	\end{figure}
	\smallbreak\noindent
	Since the limiting weights generated by the algorithm remain strictly positive, they can be used to scalarize the multi-objective problem of maximizing or minimizing all features, thereby allowing recovery of Pareto–optimal points from the dataset.
    \begin{corollary}[Pareto optimality of the equilibrium weights]\label{cor0}
		Assume that the normalization is coordinatewise order–preserving, i.e., for each 
		$j\in\{1,\dots,m\}$ and all $i,k\in\{1,\dots,n\}$
		\[
\left\{
\begin{aligned}
x_{ij} \le x_{kj} &\implies \Phi_{ij} \le \Phi_{kj},\\
x_{ij} < x_{kj} &\implies \Phi_{ij} < \Phi_{kj}.
\end{aligned}
\right.
\]
		Let $\gamma^* \in \mathcal K_+^{m-1}$ be given by (\ref{fp}) and let $\mathcal{X} = \{x_1, \dots, x_n\} \subseteq \mathbb{R}^m$ be the original dataset. Let  $F: \mathcal{X} \to \mathbb{R}^m$ be the multi-objective function given by
		\[
		F(x_i) := (x_{i1}, \dots, x_{im}).
		\]
		Then, any index $i^* \in \{1, \dots, n\}$ maximizing the scalarization
		\begin{align}\label{eq:organismfitness1}
		r_{i}:=	\sum_{j=1}^m \gamma^*_j\, \Phi_{ij}
		\end{align}
		yields a Pareto-optimal element $x_{i^*} \in \mathcal{X}$ for the problem of maximizing $F$, i.e.,
\[
F(x)\geq F(x_{i^*}) \qquad\text{implies}\qquad F(x)=F(x_{i^*})
\]
for all \(x\in\mathcal X\), where the inequality is understood componentwise. 	\end{corollary}
	The same conclusion holds if the objective is to minimize all features instead of maximize them, with the adequate changes mutatis mutandis.

	\section{Proofs}\label{proofs}

    \subsection{Proof of Lemma \ref{lem0}}
Fix $\gamma\in\mathcal K^{m-1}$ and $j\in\{1,\dots,m\}$. By definition,
\begin{align}\label{eq:sot}
1+\stepsize\Delta_j(\gamma)
&=
1-\stepsize\gamma_j\Bigl(\widetilde{\Phi}_j+\frac12\Bigr)
+\frac{2\stepsize}{m}\sum_{s=1}^m \gamma_s\widetilde{\Phi}_s.
\end{align}
Since $\sum_{s=1}^m \gamma_s\widetilde{\Phi}_s\ge \gamma_j\widetilde{\Phi}_j$ and $m\ge2$, it follows from~\eqref{eq:h0} that
\begin{align*}
1+\stepsize\Delta_j(\gamma)
&\ge
1-\stepsize\gamma_j\Bigl(\widetilde{\Phi}_j+\frac12\Bigr)
+\frac{2\stepsize}{m}\gamma_j\widetilde{\Phi}_j\ge 1-\stepsize\gamma_j\Bigl(\underbrace{\Bigl(1-\frac{2}{m}\Bigr)\widetilde{\Phi}_{\operatorname{max}}+\frac12}_{:=h_0^{-1}} \Bigr)>0.
\end{align*}
Since $\gamma\in\mathcal K^{m-1}$, there exists $r\in\{1,\dots,m\}$ such that $\gamma_r>0$. Applying the above estimate with $j=r$, we obtain  $\gamma_r\bigl(1+\stepsize\Delta_r(\gamma)\bigr)>0$. Therefore $\sum_{j=1}^m \gamma_j\bigl(1+\stepsize\Delta_j(\gamma)\bigr)>0$. This completes the proof.
\hfill$\square$
    
	\subsection{Proof of Proposition \ref{prop0}}
   Consider the log-barrier function $\Gamma:\mathcal K_+^{m-1}\to\mathbb R$ defined by
\begin{align}\label{eq:Gammadefinition}
\Gamma(\gamma):=\sum_{j=1}^m \frac{\log\gamma_j}{\widetilde{\Phi}_j+\frac12}.
\end{align}
We split the proof into four steps. In Step~1, we show that the sequence $(\gamma^k)_{k\in\mathbb N}$ remains in $\mathcal K_+^{m-1}$. In Step~2, we establish a technical inequality that will be used to control the variation of $\Gamma$ along the sequence. In Step~3, we prove that the sequence $(\Gamma(\gamma^k))_{k\in\mathbb N}$ is nondecreasing. Finally, in Step~4, we use this monotonicity to show that every accumulation point of $(\gamma^k)_{k\in\mathbb N}$ belongs to $\mathcal K_+^{m-1}$.
	\smallbreak \noindent
	\textbf{Step 1.}\,(\textit{The sequence remains in $\mathcal K_+^{m-1}$}).
We proceed by induction. By assumption, \(\gamma^0\in \mathcal K_+^{m-1}\). Assume now that \(\gamma^k\in \mathcal K_+^{m-1}\) for some \(k\in\mathbb N\). Then for every \(j\in\{1,\dots,m\}\)
we have \(\gamma_j^k>0\). We deduce from the proof of Lemma~\ref{lem0} that $
\gamma_j^k\bigl(1+h\Delta_j(\gamma^k)\bigr)>0$ and from the update rule of Algorithm~\ref{gAI} that $\gamma^{k+1}\in \mathcal K_+^{m-1}$. 
\smallbreak\noindent
\textbf{Step 2.}\,(\textit{A technical inequality}).
Fix \(\gamma\in \mathcal K_+^{m-1}\). We claim that
\begin{equation}\label{eq:tech-ineq}
\Biggl(\sum_{s=1}^m\gamma_s\bigl(1+h\Delta_s(\gamma)\bigr)\Biggr)
\sum_{j=1}^m
\frac{1}{\bigl(\widetilde{\Phi_j}+\frac12\bigr)\bigl(1+h\Delta_j(\gamma)\bigr)}
\, \le \,
\sum_{j=1}^m \left(\frac{1}{\widetilde{\Phi_j}+\frac12}\right).
\end{equation}
To prove this, we set $a_j:=\left(\widetilde{\Phi_j}+\frac12\right)^{-1}$. We need to show that  
\[
R:= \sum_{j=1}^m a_j - \Biggl(\sum_{r=1}^m\gamma_r\bigl(1+h\Delta_r(\gamma)\bigr)\Biggr)
\sum_{j=1}^m
\frac{a_j}{\bigl(1+h\Delta_j(\gamma)\bigr)} \,\ge\, 0
\]
Since
\[
1+h\Delta_j(\gamma)
=
1-h\gamma_j\Bigl(\underbrace{\widetilde{\Phi_j}+\frac12}_{:=a_j^{-1}}\Bigr)
+\frac{2h}{m}\sum_{s=1}^m\gamma_s\widetilde{\Phi_s}, \qquad \text{and}\qquad \sum_{s=1}^m
\gamma_s =1
\]
we have
\[
\sum_{s=1}^m\gamma_s\bigl(1+h\Delta_s(\gamma)\bigr)
=
1+\frac{2h}{m}\sum_{s=1}^m\gamma_s\widetilde{\Phi_r}
-
h\sum_{s=1}^m \frac{\gamma_s^2}{a_s}\,.
\]
Hence, for every \(j\in\{1,\dots,m\}\),
\begin{align}\label{eq:tri2}
1+h\Delta_j(\gamma)-\sum_{s=1}^m\gamma_s\bigl(1+h\Delta_s(\gamma)\bigr)
= h\sum_{r=1}^m
\gamma_r\Bigl(
\frac{\gamma_r}{a_r}
-
\frac{\gamma_j}{a_j}
\Bigr).
\end{align}
Writing
\[
\sum_{j=1}^m
a_j \equiv \sum_{r,j=1}^m
\frac{\gamma_r a_j (1+h\Delta_j(\gamma))}{(1+h\Delta_j(\gamma))}
\]
we eventually obtain using~\eqref{eq:tri2}:
\begin{align*}
R
&=
\sum_{r,j=1}^m
\frac{\gamma_r a_j}{\bigl(1+h\Delta_j(\gamma)\bigr)} \left(
1+h\Delta_j(\gamma)-\sum_{s=1}^m\gamma_s(1+h\Delta_s(\gamma))\right)
\\[0.5ex]
&=
\sum_{r,j=1}^m
\frac{h\,\gamma_r \,a_j}{1+h\Delta_j}\left( \Delta_j-\Delta_r\right)
= \sum_{r,j=1}^m
\frac{h\,\gamma_r\,a_j}{1+h\Delta_j}\left( \frac{\gamma_r}{a_r}-\frac{\gamma_j}{a_j}\right).
\end{align*}
Now, symmetrizing in the indices \(r\) and \(j\), we get
\begin{align*}
R
&=
\frac{h}{2} \left( 
\sum_{r,j=1}^m
\frac{\gamma_r\,a_j}{(1+h\Delta_j)}\left( \frac{\gamma_r}{a_r}-\frac{\gamma_j}{a_j}\right) + 
\sum_{r,j=1}^m \frac{\gamma_j\,a_r}{(1+h\Delta_r)} \left(\frac{\gamma_j}{a_j}-\frac{\gamma_r}{a_r}\right) \right)
\\
&=
\frac{h}{2}
\sum_{r,j=1}^m
\left( \frac{\gamma_r}{a_r}-\frac{\gamma_j}{a_j}\right)
\Biggl[
\frac{ \gamma_r\,a_j\,(1+h\Delta_r) - \gamma_j\,a_r\,(1+h\Delta_j)}
{
(1+h\Delta_j)\,(1+h\Delta_r)}
\Biggl]
\end{align*}
We set
\[
\Gamma_{r,j}:= \gamma_r\,a_j\,(1+h\Delta_r) -\gamma_j\,a_r\,(1+h\Delta_j)\,.
\]
In view of~\eqref{eq:sot} after elementary calculations we obtain
\[
\Gamma_{r,j}=\,a_ra_j\left( \frac{\gamma_r}{a_r}-\frac{\gamma_j}{a_j}\right)
\Biggl[ 
\underbrace{\left(1+\frac{2\stepsize}{m}\sum_{s=1}^m \gamma_s\widetilde{\Phi}_s\right) - h\left(\frac{\gamma_r}{a_r} +\frac{\gamma_j}{a_j} \right)}_{:=P_{r,j}}
\Biggr] 
\]
and we deduce that
\[
R=\frac{h}{2} 
\sum_{r,j=1}^m
\left( \frac{\gamma_r}{a_r}-\frac{\gamma_j}{a_j}\right)^2 \frac{a_r\,a_j\,P_{r,j}}{(1+h\Delta_r)(1+h\Delta_j)}\,.
\]
It is sufficient to prove that $P_{r,j} \ge 0$ whenever $r \neq j$. Recalling that $a_j:=\left(\widetilde{\Phi_j}+\frac12\right)^{-1}$ and noting that for \(r\neq j\) we have 
$$\sum_{s=1}^m\gamma_s\widetilde{\Phi_s}\ge \gamma_r\widetilde{\Phi_r}+\gamma_j\widetilde{\Phi_j},$$
we deduce:
\begin{align*}
&P_{r,j}=1+\frac{2h}{m}\sum_{s=1}^m\gamma_s\widetilde{\Phi_s}
-
h\bigl(\widetilde{\Phi_r}+\tfrac12\bigr)\gamma_r
-
h\bigl(\widetilde{\Phi_j}+\tfrac12\bigr)\gamma_j
\\
&\ge
1-
h\Bigl(\Bigl(1-\frac2m\Bigr)\widetilde{\Phi_r}+\frac12\Bigr)\gamma_r
-
h\Bigl(\Bigl(1-\frac2m\Bigr)\widetilde{\Phi_j}+\frac12\Bigr)\gamma_j.
\end{align*}
Since \(m\ge 2\), we have \(1-\frac2m\ge 0\), and therefore
\[
\Bigl(1-\frac2m\Bigr)\widetilde{\Phi_r}+\frac12
\le
\Bigl(1-\frac2m\Bigr)\widetilde{\Phi}_{\operatorname{max}}+\frac12,
\]
and similarly for \(j\). Hence
\begin{align*}
&1+\frac{2h}{m}\sum_{i=1}^m\gamma_i\widetilde{\Phi_i}
-
h\bigl(\widetilde{\Phi_r}+\tfrac12\bigr)\gamma_r
-
h\bigl(\widetilde{\Phi_j}+\tfrac12\bigr)\gamma_j
\\
&\ge
1-
h\left(\Bigl(1-\frac2m\Bigr)\widetilde{\Phi}_{\operatorname{max}}+\frac12\right)
(\underbrace{\gamma_r+\gamma_j}_{\le 1}).
\end{align*}
Recalling~\eqref{eq:h0}, since \(0<h<h_0\) we have $h\left(\Bigl(1-\frac2m\Bigr)\widetilde{\Phi}_{\operatorname{max}}+\frac12\right)<1$ and consequently $P_{r,j}\ge 0$.
\smallbreak\noindent
\textbf{Step 3.}\,(\textit{The sequence \((\Gamma(\gamma^k))_{k\in\mathbb N}\) is nondecreasing}).
Let \(k\in\mathbb N\). By the definition of \(\Gamma\) and the update rule of Algorithm~\ref{gAI}, we have
\begin{align*}
\Gamma(\gamma^{k+1})-\Gamma(\gamma^k)
=
\sum_{j=1}^m a_j
\log\!\left(\frac{\gamma_j^{k+1}}{\gamma_j^k}\right)
=
-A\sum_{j=1}^m \frac{a_j}{A}
\log\!\left(
\frac{\displaystyle\sum_{s=1}^m \gamma_s^k\bigl(1+h\Delta_s(\gamma^k)\bigr)}
{1+h\Delta_j(\gamma^k)}
\right)\,,
\end{align*}
for $A=a_1+\dots+a_m$. Since the function $-\log:(0,+\infty)\to\mathbb R$ is convex, we obtain
\begin{align*}
\Gamma(\gamma^{k+1})-\Gamma(\gamma^k)
\ge
-A\log\!\left(
\frac{1}{A}
\Biggl(\underbrace{\sum_{s=1}^m \gamma_s^k\bigl(1+h\Delta_s(\gamma^k)\bigr)\Biggr)
\sum_{j=1}^m
\frac{a_j}{1+h\Delta_j(\gamma^k)}}_{\leq A:=\sum_{j=1}^m a_j \,\,\text{(by Step 2)}}
\right)\,\geq 0\,.
\end{align*}
We conclude that $\Gamma(\gamma^{k+1})\ge \Gamma(\gamma^k)$ for all $k\in\mathbb N$. Thus the sequence \((\Gamma(\gamma^k))_{k\in\mathbb N}\) is nondecreasing.
\smallbreak\noindent
\textbf{Step 4.}\,(\textit{Any accumulation point belongs to \(\mathcal K_+^{m-1}\)}).
By Step~3, the sequence \((\Gamma(\gamma^k))_{k\in\mathbb N}\) is nondecreasing. Hence
\[
\Gamma(\gamma^k)\ge \Gamma(\gamma^0)
\qquad \forall k\in\mathbb N.
\]
Fix \(j\in\{1,\dots,m\}\). By Step~1, we have \(\gamma_s^k\in(0,1]\) for all \(s\in\{1,\dots,m\}\) and all \(k\in\mathbb N\). Therefore
\[
\log\gamma_s^k\le 0
\qquad \forall s\in\{1,\dots,m\},\ \forall k\in\mathbb N.
\]
Since \(a_s:=\left(\widetilde{\Phi_s}+1/2\right)^{-1}>0\) for every \(s\in\{1,\dots,m\}\), it follows that
\[
\Gamma(\gamma^k)
=
\sum_{s=1}^m a_s \log\gamma_s^k \,
\le \,
a_j \log\gamma_j^k.
\]
Combining this with \(\Gamma(\gamma^k)\ge \Gamma(\gamma^0)\), we obtain
\[
a_j\log\gamma_j^k \ge \Gamma(\gamma^0) \qquad\text{and consequently} \qquad\log\gamma_j^k \,\ge\, \frac{\Gamma(\gamma^0)}{a_j}.
\]
It follows that
\[
\gamma_j^k
\ge
\exp\!\left(\frac{\Gamma(\gamma^0)}{a_j}\right)
=:\eta_j>0
\qquad \forall k\in\mathbb N.
\]
Thus every coordinate of \(\gamma^k\) stays uniformly away from \(0\). 
\smallbreak\noindent
Let now \(\bar\gamma\in\mathcal K^{m-1}\) be any accumulation point of \((\gamma^k)_{k\in\mathbb N}\). Then there exists a subsequence \((\gamma^{k_\ell})_{\ell\in\mathbb N}\) such that
\[
\gamma^{k_\ell}\longrightarrow \bar\gamma
\quad\text{as}\quad\ell\to+\infty.
\]
Passing to the limit in the inequality \(\gamma_j^{k_\ell}\ge \eta_j\), we get
\[
\bar\gamma_j\ge \eta_j>0
\qquad \forall j\in\{1,\dots,m\}.
\]
Since also \(\sum_{j=1}^m \bar\gamma_j=1\), we conclude that \(\bar\gamma\in\mathcal K_+^{m-1}\). This proves that every accumulation point of \((\gamma^k)_{k\in\mathbb N}\) belongs to \(\mathcal K_+^{m-1}\). \hfill\(\square\)

	\subsection{Proof of Theorem \ref{thm0}}
We first prove \((a)\). Clearly,
\[
V(\gamma^*)=\sum_{j=1}^m \gamma_j^*\log\Big(\frac{\gamma_j^*}{\gamma_j^*}\Big)=0.
\]
Now let $\gamma\in\mathcal K_+^{m-1}$ and set
\[
t_j:=\frac{\gamma_j}{\gamma_j^*}>0
\qquad \forall j\in\{1,\dots,m\}.
\]
Using the elementary inequality $\log t\le t-1$ for every $t>0$, we get
\[
-\log t_j\ge 1-t_j
\qquad \forall j\in\{1,\dots,m\}.
\]
Therefore
\begin{align*}
V(\gamma)
&=
\sum_{j=1}^m \gamma_j^*\log\Big(\frac{\gamma_j^*}{\gamma_j}\Big)
=
-\sum_{j=1}^m \gamma_j^*\log t_j \\
&\ge
\sum_{j=1}^m \gamma_j^*(1-t_j)
=
\sum_{j=1}^m \gamma_j^*-\sum_{j=1}^m \gamma_j
=1-1=0.
\end{align*}
Moreover, equality holds if and only if $\log t_j=t_j-1$ for every $j$, hence if and only if
$t_j=1$ for every $j$, that is, $\gamma_j=\gamma_j^*$ for all $j$. Thus
\[
V(\gamma)>0
\qquad \forall \gamma\in\mathcal K_+^{m-1}\setminus\{\gamma^*\}.
\]
We now prove \((b)\). Set $\lambda:=1/A$. For every $\gamma\in\mathcal K_+^{m-1}$,
\begin{align*}
V(\gamma)
&=
\sum_{j=1}^m \gamma_j^*\log(\gamma_j^*)
-
\sum_{j=1}^m \gamma_j^*\log(\gamma_j)=
\sum_{j=1}^m \gamma_j^*\log(\gamma_j^*)
-
\lambda\sum_{j=1}^m a_j \log(\gamma_j)
\end{align*}
Recalling the log-barrier function $\Gamma$ defined by~\eqref{eq:Gammadefinition} (see proof of Proposition~\ref{prop0}), we obtain
\[
V(\gamma)=\sum_{j=1}^m \gamma_j^*\log(\gamma_j^*)-\lambda\,\Gamma(\gamma).
\]
In particular,
\[
V(\gamma^{k+1})-V(\gamma^k)
=
-\lambda\bigl(\Gamma(\gamma^{k+1})-\Gamma(\gamma^k)\bigr).
\]
Since $\lambda>0$ and, by Step $3$ of Proposition \ref{prop0}, the sequence $(\Gamma(\gamma^k))_{k\in\mathbb N}$ is nondecreasing, it follows that
\[
V(\gamma^{k+1})\le V(\gamma^k)
\qquad \forall k\in\mathbb N.
\]
It remains to prove the convergence $\gamma^k\to\gamma^*$. Let $T:\mathcal K_+^{m-1}\to\mathcal K_+^{m-1}$ 
be the update map of Algorithm \ref{gAI}, i.e., 
\begin{align}\label{eq:Tmapdefinition}
T_j(\gamma):=
\frac{\gamma_j(1+\stepsize\Delta_j(\gamma))}
{\displaystyle\sum_{s=1}^m\gamma_s(1+\stepsize\Delta_s(\gamma))}
\quad \forall j\in\{1,\dots,m\}.
\end{align}
By Lemma \ref{lem0}, this map is well defined, and it is continuous on $\mathcal K_+^{m-1}$. 
Since $V(\gamma^k)\ge 0$ by part \((a)\) and $(V(\gamma^k))_{k\in\mathbb N}$ is nonincreasing by part \((b)\), there exists
$\ell\ge 0$ such that
\[ 
V(\gamma^k)\longrightarrow \ell. 
\]
Because $(\gamma^k)_{k\in\mathbb N}\subset \mathcal K^{m-1}$ and $\mathcal K^{m-1}$ is compact, the sequence admits an accumulation point. Let $\bar\gamma$ be one of them. By Proposition \ref{prop0}, we have $\bar\gamma\in \mathcal K_+^{m-1}.$ 
Choose a subsequence $(\gamma^{k_r})_{r\in\mathbb N}$ such that
\[
\gamma^{k_r}\longrightarrow \bar\gamma.
\]
Since also $V(\gamma^{k_r})\to \ell$ and $V(\gamma^{k_r+1})\to \ell$, continuity of $T$ and $V$ yields
\[
0
=
\lim_{r\to\infty}\bigl(V(\gamma^{k_r+1})-V(\gamma^{k_r})\bigr)
=
V(T(\bar\gamma))-V(\bar\gamma).
\]
We claim that this implies $\bar\gamma=\gamma^*$. Indeed, let $\gamma\in\mathcal K_+^{m-1}$ be arbitrary. We set: 
$$ S(\gamma):=\sum_{s=1}^{m} \gamma_s(1+\Delta_s(\gamma)).$$ 
Then in Step~3 of the proof of Proposition~\ref{prop0}, it was shown that
\begin{align*}
\Gamma(T(\gamma))-\Gamma(\gamma)
&=
-A\sum_{j=1}^m \frac{a_j}{A}
\log\!\left(\frac{S(\gamma)}{1+\stepsize\Delta_j(\gamma)}\right)\\
&\ge
-A\log\!\left(
\frac{1}{A}S(\gamma)\sum_{j=1}^m\frac{a_j}{1+\stepsize\Delta_j(\gamma)}
\right)
\ge 0.
\end{align*}
Hence
\[
V(T(\gamma))-V(\gamma)
=
-\lambda\bigl(\Gamma(T(\gamma))-\Gamma(\gamma)\bigr)
\le 0.
\]
Assume now that
\[
V(T(\gamma))-V(\gamma)=0.
\]
Then necessarily $\Gamma(T(\gamma))-\Gamma(\gamma)=0$, so equality holds in the previous chain of inequalities. In particular, equality holds in Jensen's inequality, and since the logarithm is strictly concave, we get
\[
\frac{S(\gamma)}{1+\stepsize\Delta_1(\gamma)}
=
\cdots
=
\frac{S(\gamma)}{1+\stepsize\Delta_m(\gamma)}.
\]
Thus
\[
1+\stepsize\Delta_1(\gamma)=\cdots=1+\stepsize\Delta_m(\gamma).
\]
Recalling that
\[
\stepsize\Delta_j(\gamma)
=
-\stepsize\gamma_j\Bigl(\widetilde{\Phi}_j+\frac12\Bigr)
+\frac{2\stepsize}{m}\sum_{s=1}^m\gamma_s\widetilde{\Phi}_s,
\]
we infer that
\[
\gamma_j\Bigl(\widetilde{\Phi}_j+\frac12\Bigr)
=
\gamma_r\Bigl(\widetilde{\Phi}_r+\frac12\Bigr)
\qquad \forall j,r\in\{1,\dots,m\}.
\]
Hence there exists a constant $c>0$ such that
\[
\gamma_j=\frac{c}{\widetilde{\Phi}_j+\frac12}\equiv a_j\,c
\qquad \forall j\in\{1,\dots,m\}.
\]
Since $\sum_{j=1}^m\gamma_j=1$, we obtain $
c\sum_{j=1}^m a_j=1$
that is, 
$c=\lambda.$
Therefore
\[
\gamma_j=\lambda\,a_j=\gamma_j^*
\qquad \forall j\in\{1,\dots,m\}.
\]
We have shown that
\[
V(T(\gamma))-V(\gamma)=0
\qquad\Longrightarrow\qquad
\gamma=\gamma^*.
\]
Thus every accumulation point of $(\gamma^k)_{k\in\mathbb N}$ is equal to $\gamma^*$; whence it follows that 
\[
\gamma^k\longrightarrow \gamma^*
\quad\text{as }\quad k\to+\infty.
\]
This proves item $(c)$. 
\hfill$\square$

	\subsection{Proof of Corollary \ref{cor0}}
	Define, for each $i \in \{1,\dots,n\}$, 
	\[
	r_i := \sum_{j=1}^m \gamma^*_j \Phi_{ij},
	\] 
	and choose $i^* \in \argmax_{i \in \{1,\dots,n\}} r_i$. Suppose $x_{i^*}$ is not Pareto-optimal for maximizing the multi-objective function~$F$. Then there exists $k \in \{1,\dots,n\}$ such that 
	\[
	x_{kj} \ge x_{i^*j} \quad \text{for all } j \in \{1,\dots,m\} \quad \text{and } x_{kj_0} > x_{i^*j_0} \text{ for some } j_0 \in \{1,\dots,m\}.
	\]
	By the coordinatewise order preservation assumption, we have 
	\[
	\Phi_{kj} \ge \Phi_{i^*j} \quad \text{for all } j \in \{1,\dots,m\}, \quad \text{and } \Phi_{kj_0} > \Phi_{i^*j_0}.
	\]
	Since $\gamma^* \in \mathcal K_+^{m-1}$, we have $\gamma^*_j > 0$ for all $j \in \{1,\dots,m\}$, hence
	\[
	r_k - r_{i^*}=\sum_{j=1}^m \gamma^*_j(\Phi_{kj}-\Phi_{i^*j})>0,
	\]
	which contradicts the maximality of $r_{i^*}$. Therefore $x_{i^*}$ is Pareto-optimal for maximizing $F$.\hfill$\square$

	\section{An illustrative numerical experiment}\label{numexp}
	We consider a real-world dataset of 15 office listings for rent in Vienna, collected from the public real estate platform \texttt{immoscout24.at}. Each listing is described by four numerical features: monthly rent (in euros), office size (in square meters), number of rooms, and if they have a balcony (which is here modeled as a binary variable with values $0$ and $1$). These features are heterogeneous in scale and meaning, so direct comparison is not meaningful.
	To address this, we normalize the data into a matrix \( \Phi \in [0,1]^{15 \times 4} \) such that higher values consistently represent more desirable attributes.
	For the rent feature, where lower values are preferred, we apply a shifted inverted normalization:
	\begin{align}
		\Phi_{i1} = 1 - \frac{x_{i1} - \min_s x_{s1}}{\max_s x_{s1}}. \label{eq:shiftedinvertednormalization}
	\end{align}
	
	For the features size and rooms (columns 2 and 3), where higher values are preferred, we apply a normalization relative to their maximum values:
	\begin{align}
		\Phi_{ij} = \frac{x_{ij}}{\max_s x_{sj}} \quad \text{for } j = 2, 3. \label{eq:standardnormalization}
	\end{align}
	This choice of normalization functions aims to map the real-world dataset to a subset of $[0,1]$ in a comparable way. In particular, equation \eqref{eq:shiftedinvertednormalization} and  \eqref{eq:standardnormalization} both map to $[\min_s x_{sj}/\max_s x_{sj}, 1]$.
	\noindent
	For the last column (balcony), which is already binary and appropriately scaled, no transformation is applied.
	This results in a matrix $\Phi$ with entries in $[0,1]$, where higher values consistently correspond to more desirable properties. 
	\begin{align*}
		\Phi = \begin{bmatrix}
			0.6160 & 0.3000 & 0.2069 & 0.0000 \\
			0.8175 & 0.2891 & 0.2759 & 0.0000 \\
			0.2529 & 1.0000 & 0.4828 & 0.0000 \\
			0.4624 & 0.7152 & 0.4138 & 0.0000 \\
			0.4228 & 0.5478 & 0.4138 & 1.0000 \\
			0.5368 & 0.4761 & 0.4138 & 0.0000 \\
			1.0000 & 0.2674 & 0.1379 & 0.0000 \\
			0.6087 & 0.3804 & 0.3448 & 0.0000 \\
			0.2180 & 0.5000 & 0.5517 & 0.0000 \\
			0.5218 & 0.3457 & 0.2759 & 0.0000 \\
			0.9356 & 0.2891 & 0.2069 & 0.0000 \\
			0.1913 & 0.8326 & 1.0000 & 1.0000 \\
			0.1938 & 0.6826 & 0.4828 & 0.0000 \\
			0.1310 & 0.7283 & 0.4828 & 0.0000 \\
			0.7498 & 0.3587 & 0.2069 & 0.0000 \\
		\end{bmatrix}
	\end{align*}
    Note that for the mean values Tab.~\ref{tab:phi-means} we set the step size to $h=1$ in this example. The original (unnormalized) data is shown in Table~\ref{tab:raw-data}. Using \(\Phi\), we apply the evolutionary algorithm described in Subsection~\ref{Althm}, iterating for \(N=10\) steps. The resulting feature weight vector \(\gamma^k\) converges rapidly to the analytical limit  \(\gamma^*\) from Theorem~\ref{thm0}, as shown in Table~\ref{tab:gamma-table}. 
	\begin{table}[h!]
		\centering
		\captionsetup[table]{labelfont=normalfont,textfont=normalfont}
		\caption{Column-wise means \(\widetilde{\Phi}_j\) of the normalized matrix \(\Phi\).}
		\label{tab:phi-means}
		\begin{tabular}{lcccc}
			\toprule
			Feature & Rent & Size & Rooms & Balcony \\
			\midrule
			Mean \(\widetilde{\Phi}_j\) & 0.5106 & 0.5142 & 0.3931 & 0.1333 \\
			\bottomrule
		\end{tabular}
	\end{table}
	\smallbreak\noindent 
	The limit point in Theorem \ref{thm0} is given by
	\[
	\gamma^*_1 = 0.2117, \quad \gamma^*_2 = 0.2109, \quad \gamma^*_3 = 0.2395, \quad \gamma^*_4 = 0.3378.
	\]
	Approximately $33\%$ of the total weight is assigned to the feature 'balcony'. The reason for this behavior can be understood as 'balcony' has by far the lowest average, since only two offices have a balcony. Hence, as discussed in Section \ref{sec:minimalexample}, the fixed point formula from (\ref{fp}) gives the highest relevance to features with the lowest average. In statistically distributed data, objectives with a low average tend to have rare and large statistical outliers. The presented evolutionary approach distributes more weight to these objectives, since it is assumed that rare traits give an organism an advantage in the evolution. Indeed, in our presented example, having a balcony is exceptionally rare, making offices with one highly desirable.
	
\begin{figure}[h!]
		\centering
\includegraphics[width=0.6\textwidth]{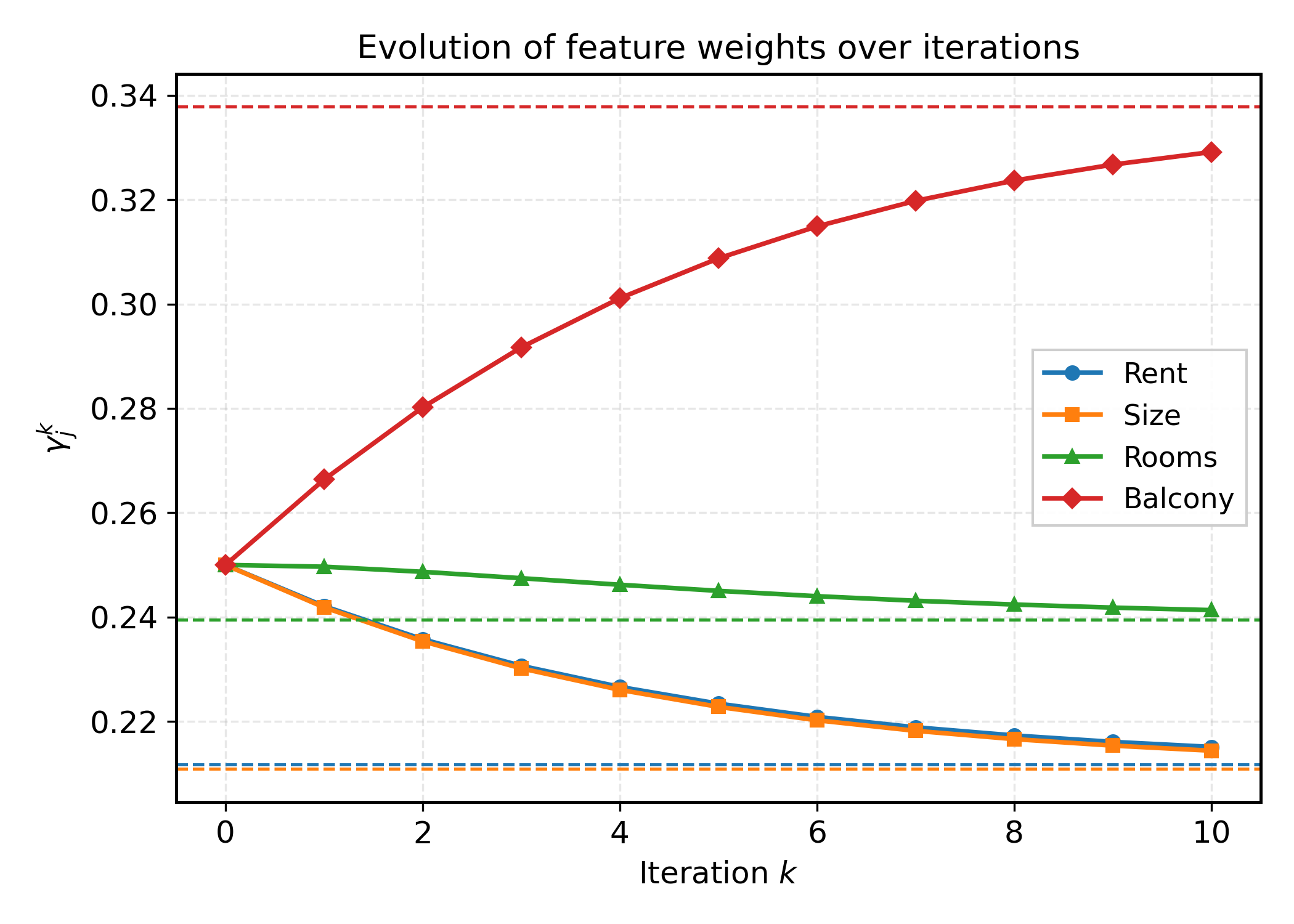}
		\caption{Evolution of feature weights $\gamma_j^k$ over $k=0,\dots,10$. 
			Dashed lines indicate the analytical solution $\gamma^*$.}
		\label{fig:gamma-evolution}
\end{figure}

	\begin{longtable}{p{6.5cm}rrrr}
		\caption{Apartment dataset including balcony.\label{tab:raw-data}} \\
		\toprule
		Name as advertised & Rent (€) & Size (m²) & Rooms & Balcony \\
		\midrule
		\endfirsthead
		
		\toprule
		Name as advertised & Rent (€) & Size (m²) & Rooms & Balcony \\
		\midrule
		\endhead
		
		\midrule
		\multicolumn{5}{r}{{Continued on next page}} \\
		\midrule
		\endfoot
		
		\bottomrule
		\endlastfoot
		
		Charmantes Altbaubüro in zentraler Lage Nähe Kärntner Straße und Oper zu mieten & 4348 & 138 & 3.0 & 0 \\
		Exklusives Büro in der Börse! & 2647 & 133 & 4.0 & 0 \\
		Wunderschöne geräumige Bürofläche im energieeffizienten Gebäude  & 7413 & 460 & 7.0 & 0 \\
		Attraktive Bürofläche in ruhiger Grünlage & 5644 & 329 & 6.0 & 0 \\
		Barrierefreie DG Bürofläche in der Favoritenstraße direkt bei der U1& 5979 & 252 & 6.0 & 1 \\
		Modernes Büro in zentraler Lage & 5016 & 219 & 6.0 & 0 \\
		Bürofläche am Höchstädtplatz: Sehr gute Infrastruktur im Haus und Umgebung | U6 & 1106 & 123 & 2.0 & 0 \\
		Repräsentatives Altbaubüro im Palais Schlick zu mieten & 4409 & 175 & 5.0 & 0 \\
		Gekühlte Bürofläche am Bauernmarkt | unbefristet & 7708 & 230 & 8.0 & 0 \\
		Gewerbefläche mit Straßenzugang in generalsaniertem Jugendstilhaus & 5143 & 159 & 4.0 & 0 \\
		Modernes Büro/Praxis in Wien: Erstbezug, 132m², U-Bahn-Nähe, Top-Ausstattung! & 1650 & 133 & 3.0 & 0 \\
		Repräsentatives Büro direkt am Karlsplatz in der Belle Etage & 7933 & 383 & 14.5 & 1 \\
		Repräsentative Bürofläche oder Praxis beim Schottentor U2/ Servitenviertel & 7912 & 314 & 7.0 & 0 \\
		Büro direkt auf der Dresdnerstraße im BC 20 - Bauteil B zu mieten & 8442 & 335 & 7.0 & 0 \\
		Bürofläche im Bürokomplex Nähe Matzleinsdorfer Platz zu mieten & 3218 & 165 & 3.0 & 0 \\
	\end{longtable}
    To quantify the deviation of the computed weights \(\gamma^*\) from the uniform baseline, we define the \emph{impact norm} as
	\[
	\|\Phi\| := \frac{m}{m-1}\sum_{j=1}^m \left| \gamma^*_j - \frac{1}{m} \right| \cdot \widetilde{\Phi}_j.
	\]
	 This measures the unevenness of weights, adjusted by the average strength of each feature. The factor $m/(m-1)$ normalizes the index so that its   maximum is $1$, attained when a single feature with average $1$ receives all the weight. 
	\smallbreak\noindent 
	To highlight the influence of the top-performing agents, we define the \emph{qualified impact norm}
	\[
	\|\Phi\|_* := \frac{m}{m-1}\sum_{j=1}^m \left| \gamma^*_j - \frac{1}{m} \right| \cdot \left( \frac{1}{|X^*|} \sum_{x_i \in X^*} \Phi_{ij} \right),
	\]
	where \(X^* \subset X\) denotes the 10\% of agents with the highest row-wise averages \(\widetilde{\phi}_i = \frac{1}{m} \sum_{j=1}^m \Phi_{ij}\). Note that this index is not a norm in the usual sense. This version emphasizes how the best listings are affected by feature weights. In this experiment, we compute $\|\Phi\| = 0.0739$ and $\|\Phi\|_* = 0.1785.$ 
	This indicates a mild deviation from uniform weighting, with more pronounced feature bias among the best-performing listings. 
	We evaluate the contribution of each feature using the \emph{feature impact} $\zeta_j := \left( \max_i \Phi_{ij} - \min_i \Phi_{ij} \right) \gamma^*_j.$
	The values for this example are given by
	\[
	\zeta_1 = 0.18397, \quad \zeta_2 = 0.15454, \quad \zeta_3 = 0.20651, \quad \zeta_4 = 0.33780.
	\]
	As for the weights $\gamma_j^*$, also for the feature impact, the feature 'balcony' dominates the example. This can be understood since an office cannot have half a balcony, but only $0$ or $1$. Hence, while other objectives show a gradual change within a certain range, 'balcony' has one single step-wise change, 'yes' or 'no'.
	\begin{table}[h!]
		\centering
		\captionsetup[table]{labelfont=normalfont,textfont=normalfont}
		\caption{Evolution of feature weights \(\gamma^k\) over $10$ iterations and comparison to the analytical solution~ \(\gamma^*\).}
        \smallbreak\smallbreak
		\label{tab:gamma-table}
		\begin{tabular}{lcccc}
			\toprule
			Iteration & Rent & Size & Rooms & Balcony \\
			\midrule
			\(\gamma^0\)   & 0.2500 & 0.2500 & 0.2500 & 0.2500 \\
			\(\gamma^1\)   & 0.2421 & 0.2419 & 0.2497 & 0.2664 \\
			\(\gamma^2\)   & 0.2358 & 0.2354 & 0.2487 & 0.2802 \\
			\(\gamma^3\)   & 0.2307 & 0.2302 & 0.2475 & 0.2917 \\
			\(\gamma^4\)   & 0.2266 & 0.2261 & 0.2462 & 0.3011 \\
			\(\gamma^5\)   & 0.2234 & 0.2228 & 0.2450 & 0.3088 \\
			\(\gamma^6\)   & 0.2209 & 0.2202 & 0.2440 & 0.3149 \\
			\(\gamma^7\)   & 0.2189 & 0.2182 & 0.2431 & 0.3198 \\
			\(\gamma^8\)   & 0.2173 & 0.2166 & 0.2424 & 0.3237 \\
			\(\gamma^9\)   & 0.2161 & 0.2154 & 0.2418 & 0.3267 \\
			\(\gamma^{10}\)& 0.2151 & 0.2144 & 0.2413 & 0.3291 \\
			\(\gamma^*\)   & 0.2117 & 0.2109 & 0.2395 & 0.3378 \\
			\bottomrule
		\end{tabular}
	\end{table}
	Figure~\ref{fig:gamma-evolution} illustrates the convergence of the weights $\gamma^k$ over the course of $10$ iterations. 
	Starting from the uniform initialization $\gamma^0=(0.25,0.25,0.25,0.25)$, each coordinate gradually drifts toward its fixed point value. 
	The weights for \emph{Rent}, \emph{Size}, and \emph{Rooms} decrease slightly and stabilize near $0.21$–$0.24$, while the weight for \emph{Balcony} steadily increases, converging to approximately $0.34$. 
	The dashed lines mark the analytical solution $\gamma^*$, showing that the iterative dynamics approach the fixed point with monotone convergence in each coordinate.
	\noindent
	For each office $x_i$, we define an aggregate score $r_i = \sum_{j=1}^m \gamma_j \,\Phi_{ij}$, 
	which evaluates the overall desirability of $x_i$ given a weight vector $\gamma$.  
	With uniform weights, $\gamma^{\operatorname{uniform}}=(0.25,0.25,0.25,0.25)$, all features contribute equally and produce the baseline scores $r_i^{\text{init}}$.  
	With $\gamma=\gamma^*$, the fixed point solution, the scores $r_i^{\text{evol}}$ incorporate the learned feature importances.  
	Ranking the offices by $r_i$ under these two regimes gives the comparison shown in Table~\ref{tab:ranking}.
	\begin{table}[H]
		\centering
		\footnotesize
		\captionsetup[table]{labelfont=normalfont,textfont=normalfont}
		\caption{Office rankings under uniform weights $\gamma^{\mathrm{uniform}}=(0.25,0.25,0.25,0.25)$ and under $\gamma^*=(0.2117,0.2109,0.2395,0.3378)$.}
        \smallbreak\smallbreak
		\label{tab:ranking}
		\begin{tabular}{@{}r|r|rrrr|r|rrrr@{}}
			\toprule
			& \multicolumn{5}{c|}{$\gamma^{\mathrm{uniform}}$} & \multicolumn{5}{c}{$\gamma^*$} \\
			Rank & Score & Rent & Size & Rooms & Balc. & Score & Rent & Size & Rooms & Balc. \\
			\midrule
			1  & 0.755979 & 7933 & 383 & 14.5 & 1 & 0.793484 & 7933 & 383 & 14.5 & 1 \\
			2  & 0.596097 & 5979 & 252 & 6 & 1 & 0.641988 & 5979 & 252 & 6 & 1 \\
			3  & 0.433915 & 7413 & 460 & 7 & 0 & 0.380130 & 7413 & 460 & 7 & 0 \\
			4  & 0.397865 & 5644 & 329 & 6 & 0 & 0.347898 & 5644 & 329 & 6 & 0 \\
			5  & 0.357897 & 1650 & 133 & 3 & 0 & 0.313203 & 5016 & 219 & 6 & 0 \\
			6  & 0.356680 & 5016 & 219 & 6 & 0 & 0.308616 & 1650 & 133 & 3 & 0 \\
			7  & 0.351331 & 1106 & 123 & 2 & 0 & 0.301152 & 1106 & 123 & 2 & 0 \\
			8  & 0.345613 & 2647 & 133 & 4 & 0 & 0.300664 & 7912 & 314 & 7 & 0 \\
			9  & 0.339790 & 7912 & 314 & 7 & 0 & 0.300134 & 2647 & 133 & 4 & 0 \\
			10  & 0.335508 & 8442 & 335 & 7 & 0 & 0.297003 & 8442 & 335 & 7 & 0 \\
			11  & 0.333501 & 4409 & 175 & 5 & 0 & 0.291728 & 4409 & 175 & 5 & 0 \\
			12  & 0.328854 & 3218 & 165 & 3 & 0 & 0.283968 & 3218 & 165 & 3 & 0 \\
			13  & 0.317420 & 7708 & 230 & 8 & 0 & 0.283780 & 7708 & 230 & 8 & 0 \\
			14  & 0.285828 & 5143 & 159 & 4 & 0 & 0.249463 & 5143 & 159 & 4 & 0 \\
			15  & 0.280716 & 4348 & 138 & 3 & 0 & 0.243249 & 4348 & 138 & 3 & 0 \\
			\bottomrule
		\end{tabular}
	\end{table}

	\section{Analogies with genetic algorithms and evolutionary systems}
	\label{sec:Analogies}
	Given the fixed point \(\gamma^*\) from Theorem~\ref{thm0}, a natural question is how the update functions \(\Delta_j^{\operatorname{dom}}\), \(\Delta_j^{\operatorname{bal}}\), and the overall update rule are derived. In what follows, we provide an evolutionary interpretation of the matrix \(\Phi\), which leads directly to the definitions in Section~\ref{Althm}.
	
	\smallbreak
	
	We interpret each agent \(x_i \in \mathcal{X}\) as an organism in an evolutionary environment. The goal of the optimization is then understood as identifying the fittest organism. A natural definition for the \emph{global fitness} of organism \(i\) is
	\begin{align}
		r_i := \sum_{j=1}^m \gamma_j \Phi_{ij}, \label{eq:organismfitness}
	\end{align}
	where \(\gamma \in \mathcal{K}^{m-1}\) is a weight vector assigning relevance to each feature. We assume that the normalized matrix \(\Phi\) is constructed so that higher values of \(\Phi_{ij}\) reflect better performance of organism~\(i\) in feature \(j\). As all entries of \(\Phi\) lie in \([0,1]\), it follows that \(r_i \in [0,1]\) as well. 
	Under this interpretation, the entry \(\Phi_{ij}\) represents the \emph{local fitness} of organism \(i\) with respect to feature \(j\). From an evolutionary perspective, we may view each feature as a gene, and \(\gamma_j\) as the fitness or strength of expression of gene \(j\). The global fitness \(r_i\) thus reflects how well the organism performs, taking into account both the quality of its features and the importance assigned to them by \(\gamma\).
	\smallbreak

	In Dawkins~\cite{daw_2006}, one central argument states that genes and not organisms are the 'unit' of competition in a species. That means that the more the phenomenological features of a gene help the overall fitness of an organism, the more successful a gene will be in spreading to a population through replication. In newer editions of~\cite{daw_2006} and also in the follow-up Dawkins~\cite{daw_1982}, this focus on genes is somehow reduced and the trade-off between gene and organism behavior comes stronger into play.
	
	One of the central paradigms of biological entities in evolution is that they have to efficiently use the available resources. From an organism point of view,  this means that there has to be a trade-off between beneficial features for reproduction and their cost in terms of resources or energy consumption. For example, if the energy cost to grow wings is too substantial compared to the gain of collected food, evolution will not favor their development over the course of generations. From a mathematical perspective, this process is nothing more than an optimization of genes and organisms with a fixed overall amount of resources.
	
	In the spirit of the work of Dawkins and in our evolutionary picture of data, let us conduct a 'thought experiment': what if, to analyze the fitness of genes and organisms, we convert $\Phi$ into an evolutionary simulation, where we iteratively test local and global fitness? To that end, we envision a virtual, dimensionless 'heap of resources' and calculate the share of resources a gene obtains. The more shares a gene can accumulate, the more fit it will be in the next iteration, i.e., the larger the fitness of the gene $\gamma_j$.
	
	Continuing with an absolute heap size, would lead to unstable and inconsistent behavior, since the order genes accessing the heap would influence the result. Instead, we continue with a differential share $\Delta_j$ for the $j$th gene and require (i) larger positive shares $\Delta_j$ tend to increase $\gamma_j$ and vise versa, (ii) the total sum of gene fitness values is constant, (iii) the gene fitness does not change too abruptly to allow a fair comparison between genes. Taking these conditions into account, we can formulate a genetic replicator equation~\cite{MaynardSmith1982, Hofbauer1998}:
	\begin{align}
		\gamma_j^{k+1}=\frac{\gamma_j^k \left( 1 + \Delta_j^k \right)}{\displaystyle\sum_{s=1}^m \gamma_s^k \left( 1 + \Delta_s^k\right)} \quad \forall j\in \{1, \dots, m\}. \label{eq:replicatorequation}
	\end{align}
	
	As a next step in our evolutionary representation, we have to define the functions $\Delta_j^k$. As described above, in biology, there is an interplay between the genes and organisms in evolution. On the one hand, evolution in connection with genes can be interpreted as \emph{local} comparison of fitness in terms of a single feature. On the other hand, the evolution of organisms in a species can be understood as a \emph{global} comparison with respect to all genes available in the population. Hence, also in our evolutionary interpretation there are two contributions to $\Delta_j$
	\begin{align*}
		\Delta_j^k = \Delta_j^{k,\operatorname*{gen}} + \Delta_j^{k,\operatorname*{org}}
	\end{align*}
	where $\Delta_j^{k,\operatorname*{gen}}$ denotes the gene contribution and $\Delta_j^{k,\operatorname*{org}}$ denotes the organism contribution, respectively.  
	
	The functions $\Delta_j^{k,\operatorname*{gen}}$ and $\Delta_j^{k,\operatorname*{org}}$ govern the dynamic behavior of the system via the replicator equations~\eqref{eq:replicatorequation}. In the evolutionary picture these functions correspond to strategies of individual genes and organisms. Hence, the function are assembled as
	\begin{align}
		\Delta_j^{k,\operatorname*{gen}}&=\frac{1}{n}\sum_{i=1}^n \Delta_{ij}^{k,\operatorname*{gen}},\label{eq:avereragegenedelta}\\
		\Delta_j^{k,\operatorname*{org}}&= \frac{1}{n}\sum_{i=1}^n \Delta_{ij}^{k,\operatorname*{org}}\label{eq:avererageorgdelta},
	\end{align}
	where $\Delta_{ij}^{k,\operatorname*{gen}}$, $\Delta_{ij}^{k,\operatorname*{org}}$ denote the contribution to the differential share of a single gene or organism, respectively, and recalling that $n$ denotes the total number of organisms in the system.
	
	In this paper, we will derive a simple variant for the strategy functions $\Delta_{ij}^{k,\operatorname*{gen}}$ and $\Delta_{ij}^{k,\operatorname*{org}}$, leaving more complex strategies for future work. Generally, as conditions for the choice of strategies, we require that not all gene fitness be shifted to a single gene. In that case a single $\gamma_{j'}=1$ while all other weights in equation \eqref{eq:organismfitness} would be zero. In the optimization this would not be desirable, since the information of $\Phi_{ij}$ for $j\neq j'$ would be lost.
	
	\subsection{Gene contribution $\Delta_{ij}^{k,\operatorname*{gen}}$}
	Consider a gene describing a single feature in a real evolutionary system, e.g., the size of an animal in a species. Further, assuming that larger size is beneficial in the chosen ecosystem, the fitness of the gene 'size' in a population can be estimated by taking into account the sizes of all animals in the population and comparing them to genes describing other features of the animal. Qualitatively speaking, if for example many large animals in a population are to be found, one can deduct that size is genetically more important than other features therein.
	
	In addition to the magnitude of a quantity such as size, the importance of the feature in a previous reproduction cycle also influences the future fitness of a gene in the evolution. Translating this dynamics to our optimization problem, we define '$\operatorname*{dom}$' the dominant gene strategy function\footnote{For the evolutionary interpretation we set the step size $h=1$ in the following. A global scaling factor $h$ changes the available resources for all genes and organisms equally and not the relative shares each gene obtains.}
	\begin{align}
		\Delta_{ij}^{k,\operatorname*{gen}}=\Delta_{ij}^{k,\operatorname*{dom}}:=\gamma_j^k \left(\Phi_{ij}-\frac{1}{2}\right). \label{eq:dominantstrategy}
	\end{align}
	Note that since $\Delta_{ij}^{k,\operatorname*{dom}}$ is used to compute $\gamma_j^{k+1}$ via equation \eqref{eq:replicatorequation}, the linear factor $\gamma_j^k$ in equation \eqref{eq:dominantstrategy} is the only sequence-dependence factor that takes the current gene fitness into account for the next iteration. Additionally, the factor $(\Phi_{ij}-\frac{1}{2})$ will yield a positive differential share when $\Phi_{ij}>1/2$ and vice versa.
	
	Taking only $\Delta_{ij}^{k,\operatorname*{dom}}$ for $\Delta_{j}^k$ into account, it can be easily shown that the objective $j'$ with maximum $\widetilde{\Phi_{j'}}$ will accumulate all the weight, i.e. $\gamma_{j'}=1$. Hence, for a non-trivial result, we need $\Delta_{ij}^{k,\operatorname*{org}}$ to counter the direct dominant effect of the gene strategy function $\Delta_{ij}^{k,\operatorname*{dom}}$.
	
	\subsection{Organisms contribution $\Delta_{ij}^{k,\operatorname*{org}}$}
	In evolution, organisms with only a few excellent features have difficulties to thrive. It is usually the lifeforms flexible enough to quickly adapt to a changing habitat to win the day. Applying this paradigm to our optimization problem $\Phi$, we need a quantity to measure the dependence of an organism on a particular gene
	\begin{align}
		\mu_{ij}:=\frac{\gamma_j \Phi_{ij}}{\sum_s \gamma_s\Phi_{is}}=\frac{\gamma_j \Phi_{ij}}{r_{i}}\label{eq:mu}
	\end{align}
	Note that the larger $\mu_{ij}$, the more the fitness of an organism $i$ depends on a particular gene $j$. Furthermore, from an organism perspective, with $m$ features at hand, the ideal share would be $\mu_{ij}=\frac{1}{m}$ for all $j\in{1,\ldots,m}$.
	
	Consequently, we obtain as balanced organism strategy function
	\begin{align}
		\Delta_{ij}^{k,\operatorname*{org}}=\Delta_{ij}^{k,\operatorname*{bal}}:=-2 r_i\left(\mu_{ij}-\frac{1}{m}\right). \label{eq:balancedstrategy}
	\end{align}
	Note that a perfectly balanced organism, where $\mu_{ij}=\frac{1}{m}$ for all $j\in \{1,\ldots,m\}$, will yield no differential change to any $\Delta_j^{bal}$.
	
	\subsection{Summary of evolutionary argument}
	
	It is important to observe that equation \eqref{eq:dominantstrategy} for $\Delta_{ij}^{k,\operatorname*{dom}}$ yields a \emph{relative} change to the gene fitness $\gamma_j$. In particular, for a specific organism $i$, a gene $j'$ can only gain fitness, if $\Phi_{ij'}$ is larger than the other $\Phi_{ij}$, $j\neq j'$. Thus, the effect of $\Delta_{ij}^{k,\operatorname*{gen}}$ is directly proportional to the rate of asymmetry of an agents data. In contrast, the minus sign in equation \eqref{eq:balancedstrategy} for $\Delta_{ij}^{k,\operatorname*{bal}}$ indicates that asymmetry in an agent is penalized. The total differential share $\Delta_{ij}^{k} := \Delta_{ij}^{k,\operatorname*{dom}} + \Delta_{ij}^{k,\operatorname*{bal}}$ is thus an evolutionary trade-off between these two effects related to the asymmetry of a feature $j'$ compared to other objectives in optimization.
	
	Combining all contributions for the dominant strategy, we obtain
	\begin{align*}
		\Delta_j^{k,\operatorname*{dom}}:=\frac{1}{n}\sum_{i=1}^n \Delta_{ij}^{k,\operatorname*{dom}} = \frac{1}{n}\sum_{i=1}^n \gamma_j\left(\Phi_{ij}-\frac{1}{2}\right)=\gamma_j\left(\widetilde{\Phi_{j}}-\frac{1}{2}\right)
	\end{align*}
	with the mean value $\widetilde{\Phi_{j}}$, leading into the corresponding definition of Sec.~\ref{Althm}. Analogously, we accumulate all contributions for the balanced strategy as 
	
	\begin{align*}
		\Delta_j^{k,\operatorname*{bal}}:=\frac{1}{n}\sum_{i=1}^n \Delta_{ij}^{k,\operatorname*{bal}} = - \frac{2}{n}\sum_{i=1}^n r_i\left(\mu_{ij}-\frac{1}{m}\right)=- \frac{2}{n}\sum_{i=1}^n \left(\gamma_j \Phi_{ij}-\frac{1}{m}\sum_s \gamma_s\Phi_{is}\right)
	\end{align*}
	again leading to the corresponding expression in Sec.~\ref{Althm}. 
	
	\subsection{Asymptotic behavior in higher dimensions}
	
	The update terms \(\Delta_j^{\operatorname{dom}}\) and \(\Delta_j^{\operatorname{bal}}\) can be interpreted as two competing evolutionary pressures: the former rewards dominance of feature \(j\) across the population, while the latter penalizes imbalance in how much feature \(j\) contributes to the overall fitness of agents. For the algorithm to remain stable and unbiased as the number of features \(m\) grows, both terms should scale similarly with respect to \(m\).
	To that end, we examine the asymptotic behavior of the two update components.
	\smallbreak 
	As the number of features \(m\) increases, and in the absence of strong preference for any particular one, the weights \(\gamma_j\) are expected to distribute more evenly across all features. This implies that for large \(m\), each \(\gamma_j\) is typically of order \(\mathcal{O}(1/m)\). 
	Recall that the dominance term is given by
	\[
	\Delta_j^{\operatorname{dom}} = \gamma_j \cdot \left(\widetilde{\Phi_j} - \frac{1}{2} \right),
	\]
	where \(\widetilde{\Phi_j}\) is the average value of feature \(j\) across all agents. The quantity \(\widetilde{\Phi_j} - \tfrac{1}{2}\) measures how much the average of feature \(j\) deviates from a neutral baseline. Since \(\widetilde{\Phi_j}\) remains in the compact set $[0,1]$ as \(m\) increases, the product \(\gamma_j(\widetilde{\Phi_j} - \tfrac{1}{2})\) inherits the scaling \(\mathcal{O}(1/m)\) from \(\gamma_j\). We therefore conclude that
	\[
	\Delta_j^{\operatorname{dom}} = \mathcal{O}\left(\tfrac{1}{m}\right).
	\]
	This scaling reflects that, as the number of features grows, the influence of each individual feature on the evolutionary update diminishes—unless that feature significantly deviates from the average.
	The balancing term is defined as
	\[
	\Delta_j^{\operatorname{bal}} = -\beta\left(\gamma_j\widetilde{\Phi_j} - \frac{1}{m}\sum_{s=1}^m \gamma_s \widetilde{\Phi_s} \right),
	\]
	where the factor \(\beta = 2\) is chosen to ensure that this penalty has the same scale as the dominance term. The term inside the parentheses measures how much the contribution of feature \(j\) deviates from the average over all features. Since the mean \(\sum_{s=1}^m \gamma_s \widetilde{\Phi_s}\) is of order \(\mathcal{O}(1)\), both terms in the difference are \(\mathcal{O}(1/m)\), so the result is again
	\[
	\Delta_j^{\operatorname{bal}} = \mathcal{O}\left(\tfrac{1}{m}\right).
	\]
	The choice \(\beta = 2\) reflects a symmetry: since \(\Delta_j^{\operatorname{dom}}\) is centered around \(1/2\), the balancing correction must be scaled accordingly to counteract overdominance and ensure convergence to a stable equilibrium. More generally, if the dominance term were centered around a constant \(a\), one would require \(\beta = 1/a\) to retain this matching.
	Note that the introduced dimensional argument is a good analysis tool also for more complex evolutionary strategies. In general, all strategies included in an algorithm should show the same asymptotic behavior in terms of $m$ to allow for a fair comparison.

	\subsection{A minimal example and its evolutionary interpretation}\label{sec:minimalexample}
	
	Let us investigate the effect of the evolutionary approach on a minimal $2\times2$ example
	\begin{align*}
		X:= \begin{bmatrix}
			1&0\\
			0.5&0.5
		\end{bmatrix}.
	\end{align*}
	Since the entries are already in the interval $[0,1]$, we take $\Phi = X$. 
	Inserting the mean values $\widetilde{\Phi_1}=0.75$ and $\widetilde{\Phi_2}=0.25$ into the fixed point formula \eqref{fp} yields
	\begin{align*}
		\gamma_1^*=0.375, \quad\gamma_2^*=0.625
	\end{align*}
	Most of the relevance is assigned to the second feature, which has the lower average. In our evolutionary interpretation, this means that a lower average local fitness for a given gene makes its rare, larger values even more significant.
	\smallbreak \noindent 
	As before, imagine an animal population where size is a relevant feature described by a specific gene. Furthermore, assume that the majority of animals in the population are small. Then, a single larger organism has a decisive evolutionary advantage, since there are only a few competitive opponents in terms of size when fighting for the available resources.
	\smallbreak \noindent
	Let us now investigate variants of the previous example. For each $\xi\in[0,0.5]$, consider the $2\times 2$ matrix 
	\begin{align*}
		X(\xi):= \begin{bmatrix}
			1&0\\
			0.5+\xi&0.5-\xi 
		\end{bmatrix},
	\end{align*}
	As before, take $\Phi(\xi) = X(\xi)$. Note that as $\xi\rightarrow 0.5$ the matrix takes its maximally  asymmetric form
	\begin{align*}
		X(0.5)= \begin{bmatrix}
			1&0\\
			1&0
		\end{bmatrix}.
	\end{align*}
	Table~\ref{tab:minimalexample} shows the values of $\gamma^*$ as a function of $\xi$, showing linear behavior with respect to $\xi$ with slope  $0.25$.
	\begin{figure}[h!]
		\centering
		\includegraphics[width=0.6\textwidth]{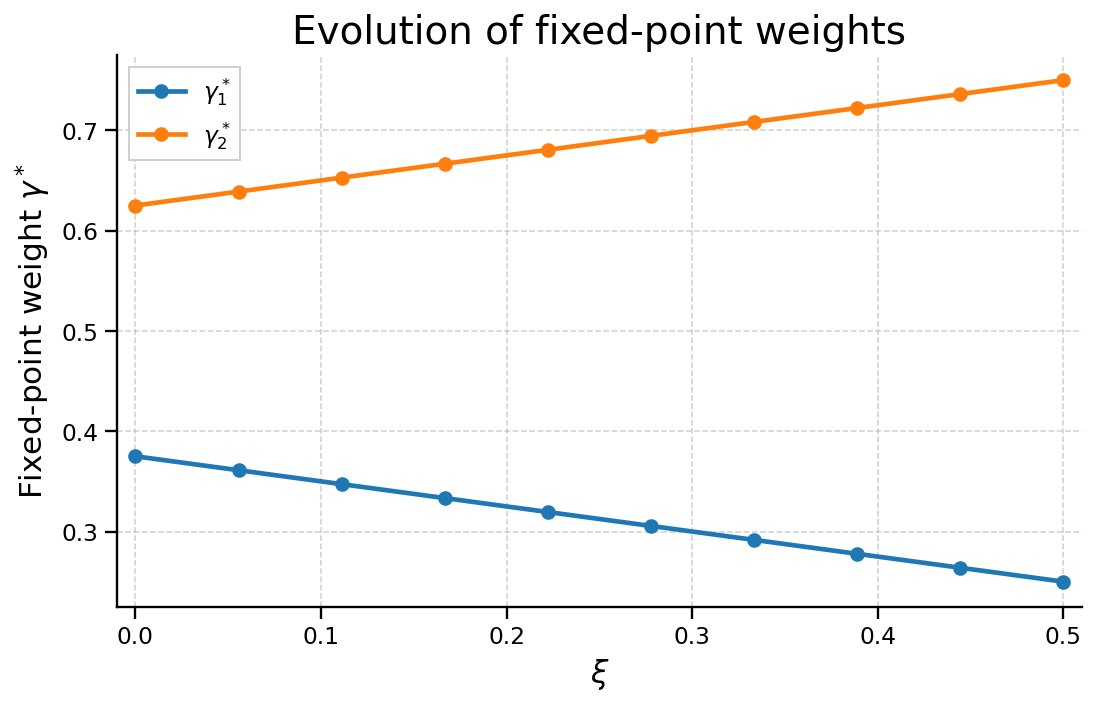}
		\caption{Evolution of fixed-point weights $\gamma_1^*$ and $\gamma_2^*$ as a function of $\xi$.}
	\end{figure}
	
	\begin{table}[h!]
    \centering
    \captionsetup[table]{labelfont=normalfont,textfont=normalfont}
    \caption{Evolution of feature weights $\gamma$ for $6$ values of $\xi$.}
    \label{tab:minimalexample}
    \setlength{\tabcolsep}{12pt} 
    \begin{tabular}{lcc}
        \toprule
        $\xi$ & $\gamma_1^*$ & $\gamma_2^*$ \\
        \midrule
        0   & 0.375 & 0.625 \\
        0.1 & 0.350 & 0.650 \\
        0.2 & 0.325 & 0.675 \\
        0.3 & 0.300 & 0.700 \\
        0.4 & 0.275 & 0.725 \\
        0.5 & 0.250 & 0.750 \\
        \bottomrule
    \end{tabular}
\end{table}
	The final values $\gamma_1^*=0.25$ and $\gamma_2^*=0.75$ are the maximally asymmetric weights for the two features.
\begin{remark}
    In the extreme case $\xi = 0.5$, both features are constant across all samples. 
    From a purely discriminative perspective, such features carry no information and 
    should be assigned equal (or zero) weight. However, the evolutionary weighting 
    rule still assigns $\gamma_1^* = 0.25$ and $\gamma_2^* = 0.75$, reflecting its 
    inherent bias toward features with lower means. One may view the situation 
    shortly before everything collapses to $1$ or $0$ as the more meaningful case, 
    where the rule amplifies small deviations as rare traits. In the trivial case the unequal weights can be seen as a consequence of continuity of the formula. 
\end{remark}

    \subsection{Relation of evolutionary picture and input data}
    The method introduced in this paper, relates the raw input data $X$ with a well-defined equilibrium~$\gamma^*$. 
    The process can be summarized as
    \begin{align*}
    \text{input data } X \rightarrow \text{normalized data } \Phi \rightarrow \text{evolutionary fitness} \rightarrow \gamma^*.
    \end{align*}
   Although some information is lost by using only the mean values $\widetilde{\Phi_{j}}$ to compute $\gamma^*$, there is often statistical value in comparing different objectives via $\widetilde{\Phi_{j}} \leftrightarrow \widetilde{\Phi_{j'}}$. In particular, differences in $\widetilde{\Phi_{j}}$ for $j \in \{1,\ldots,m\}$ often indicate distinct clustering of features in the original data $X$. 
    In this sense, the introduced evolutionary approach provides a means to assess the relevance of features in multi-objective optimization. It is important to emphasize that this evaluation is not rule-based but depends solely on the input data $X$ and the dynamics of the associated evolutionary representation.
	
	\subsection{Related literature}\label{rellit}
	Our algorithm determines (learns) a probability vector of feature weights by a multiplicative (replicator–type) update on the simplex and admits a closed-form interior limit. This places the algorithm between evolutionary multi-objective optimization, which typically evolves solutions, and feature-weighting methods in machine learning, which estimate weights but rarely via a dynamical system. 
	\smallbreak \noindent 
   We refer to the book \cite{MaynardSmith1982} for background on evolutionary game theory, and now turn to the related literature on evolutionary multi-objective algorithms. In evolutionary multi-objective optimization (EMO), algorithms mainly differ in the way they approximate the Pareto front. A comparison of several such methods is given in \cite{zitzler2000comparison}, where common design principles are identified, including Pareto-dominance ranking, external archives, performance indicators, and decomposition into scalar subproblems. Algorithms such as the Non-dominated Sorting Genetic Algorithm II (NSGA-II) \cite{deb2002nsga2} and the Strength Pareto Evolutionary Algorithm 2 (SPEA2) \cite{zitzler2001spea2} rely on Pareto-based strategies without incorporating weight vectors. In contrast, the Multi-Objective Evolutionary Algorithm based on Decomposition (MOEA/D) \cite{zhang2007moead} explicitly uses weight vectors to decompose a multi-objective problem into scalar subproblems. These weight vectors determine search directions and, to a large extent, the distribution of the final solution set \cite{ma2020survey}. In its original form, MOEA/D uses fixed weight vectors, but subsequent research has introduced adaptive strategies.
	For example, pa\(\lambda\)-MOEA/D (Pareto-adaptive MOEA/D) \cite{jiang2011pal} adjusts weight vectors using geometric features of the estimated Pareto front; DMOEA/D (Diversity-maintained MOEA/D) \cite{gu2012dmoead} employs projection/equidistant interpolation of reference (weight) vectors from the current nondominated set to preserve diversity; and MOEA/D with Adaptive Weight Adjustment (MOEA/D-AWA) \cite{qi2014moeadawa} periodically removes weight vectors in crowded regions and inserts new ones in sparse regions during the search.
	Other notable contributions include the Reference Vector Guided Evolutionary Algorithm (RVEA) \cite{cheng2016rvea} and its later variants. 
	Our approach diverges from the ones previously mentioned. Instead of adapting weight vectors to shape search directions, we evolve feature weights from data, with the interior equilibrium itself providing a scalarizer.
	\smallbreak\noindent
	In machine learning, classical filter methods assign fixed importances to features. For instance, Relief \cite{kira1992relief} uses nearest‑neighbor contrasts, and ReliefF \cite{kononenko1994relieff} extends this to multi‑class and noisy data. Another popular method, mRMR (short for Minimum Redundancy – Maximum Relevance), selects features that are both highly informative about the target and minimally redundant with each other, typically using mutual information or F‑statistics for relevance and pairwise redundancy measures \cite{peng2005mrmr}. In unsupervised learning, feature-weighted $k$-means methods span decades of variants \cite{deamorim2016survey}. Evolutionary strategies are also present, e.g., estimation-of-distribution methods applied to feature weighting \cite{inza2000fsseba}, and comparisons between genetic and co-evolutionary schemes \cite{blansche2005micai}. These methods generally optimize weights empirically but do not define a replicator-style dynamical system nor provide closed-form equilibria as in our approach. In \cite{Demirovic2020}, representative subsets of the bi-objective Pareto front are studied. By contrast, we do not select Pareto representatives, but learn a data-dependent scalarization vector on the simplex that induces a ranking of the alternatives.
	\smallbreak\noindent
  We also briefly comment on the recent paper \cite{Li2026}. Our analysis is restricted to structured datasets \(X\), whose features are either readily interpretable or sampled from well-specified probability distributions. In this setting, a direct comparison with large-scale machine learning datasets such as the image data in \cite{Li2026} is of limited relevance, since a feature-local notion of fitness is typically absent there.

    \appendix

\section{A larger example and benchmarks}
\label{sec:largerbenchmarks}

\begin{figure}[h!]
	\centering
	\includegraphics[width=0.7\textwidth]{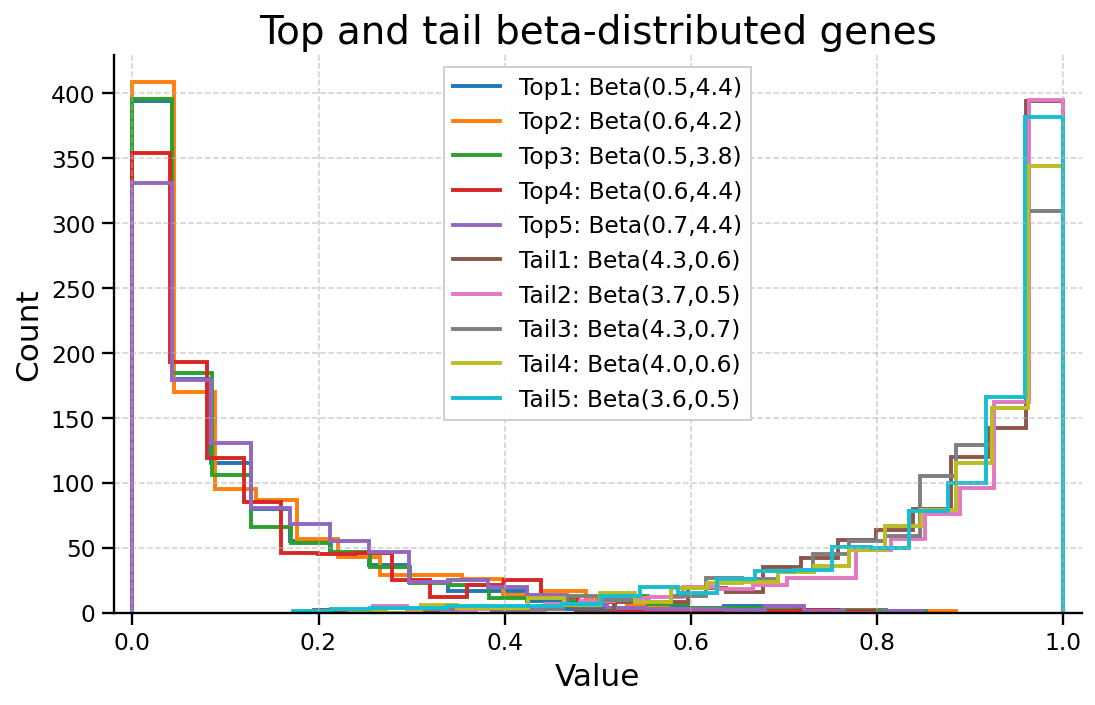}
	\caption{Histograms of the five genes with the highest and the five genes with the lowest fitness values \(\gamma_j\).}
	\label{fig:betas}
\end{figure}

We have illustrated the proposed method on small introductory examples. We now show that it also applies to larger synthetic datasets. To this end, we consider the case \(n=m=1000\), so that the system contains \(1000\) organisms and \(1000\) genes. To generate the data, we independently choose parameters \(\alpha_j,\beta_j\in(0.5,5)\) for each column \(j\in\{1,\dots,m\}\). We then draw the entries \(x_{ij}\) independently from the beta distribution \(\mathrm{Beta}(\alpha_j,\beta_j)\). After constructing the matrix \(X\), we normalize it columnwise by setting
\[
\Phi_{ij}=\frac{X_{ij}}{\max_{1\le r\le n} X_{rj}}.
\]
Applying Algorithm~\ref{gAI} to the normalized matrix \(\Phi\) with stopping tolerance \(\epsilon:=10^{-6}\), the iteration converges after \(1815\) steps. Figure~\ref{fig:betas} displays the empirical distributions of the entries \(\Phi_{ij}\) for the five genes with the largest and the five genes with the smallest limiting fitness values \(\gamma_j\).
\smallbreak\noindent
The genes with highest fitness typically correspond to columns whose beta distributions have small values of \(\alpha_j\) and relatively large values of \(\beta_j\). Such distributions are concentrated near the lower end of the interval \([0,1]\), and therefore produce smaller column averages \(\widetilde{\Phi}_j\). By the fixed-point formula \eqref{fp}, smaller averages lead to larger limiting weights \(\gamma_j\). Conversely, genes corresponding to parameters \(\alpha_j>3.5\) and \(\beta_j<1\) are associated with distributions concentrated near \(1\), hence with larger averages \(\widetilde{\Phi}_j\), and therefore with smaller limiting weights.
The larger matrix \(\Phi\) also allows us to benchmark the iterative method against the closed-form expression \eqref{fp}. Over \(1000\) independent runs of a Python implementation of Algorithm~\ref{gAI}, we obtain the timings reported in Table~\ref{tab:benchmark}.
\begin{table}[bh]
		\centering
		\captionsetup[table]{labelfont=normalfont,textfont=normalfont}
		\caption{Benchmarks comparing iterative with fixed point approach}
		\label{tab:benchmark}
		\begin{tabular}{lccc}
			\toprule
			&min[ms] & mean[ms] & max[ms]\\
			\midrule
			Iterations & 16.443  & 20.824 & 95.633 \\
            Fixed point  & \phantom{0} 0.193  &  \phantom{0}0.219 & \phantom{0}1.026\\
			\bottomrule
		\end{tabular}
	\end{table}
\noindent
For our example with $n\times m= 1000\times1000$ we find that using the fixed point formula is around $100$ times faster than iteratively solving the problem.

	

\end{document}